\documentclass[lang = american]{ems-icm} 

\theoremstyle{plain}
\newtheorem{theorem}{Theorem}

\begin{document}

\volume{?} 

\title{The work of James Maynard}
\titlemark{The work of James Maynard}

\emsauthor{1}{Kannan Soundararajan}{K.~Soundararajan}


\emsaffil{1}{Department of Mathematics, Stanford University, Stanford CA 94305 \email{ksound@stanford.edu}}


\classification[11N32, 11N35, 11J83]{11N05}

\keywords{distribution of primes, sieve methods, metric Diophantine approximation}

\begin{abstract}
 We give a brief account of some of the most spectacular results established by James Maynard, for which he has been awarded the Fields Medal.  
 \end{abstract}

\maketitle


 James Maynard has established several spectacular results in analytic number theory. While the proofs of these results involve many deep ideas, their statements are remarkable for their simplicity and elegance.   To illustrate, we state two such striking results of Maynard concerning prime numbers, before setting them in context. 

\begin{theorem} (Maynard \cite{M1}) \label{thm1} For each natural number $m \ge 2$, there exists a positive integer $C(m)$ with the following property:  There are infinitely many natural numbers $n$ such that the interval $[n, n+C(m)]$ contains at least $m$ prime numbers. 
\end{theorem}

\begin{theorem} (Maynard \cite{M2})  \label{thm2} There are infinitely many prime numbers $p$ whose decimal representation does not contain the digit $7$.
\end{theorem} 

\noindent \textbf{Background.}  To place these results in context, recall that the prime number theorem gives an asymptotic for $\pi(x)$, the number of primes below $x$; namely, 
$$ 
\pi(x) \sim \text{li}(x) = \int_0^x \frac{dt}{\log t}. 
$$
We may think of this asymptotic as roughly saying that the ``chance" of a number $n$ being prime is about $1/\log n$.  One overarching theme in analytic number theory may be formulated as asking in what ways does the sequence of primes resemble, or differ from, a random sequence of integers with each integer $n\ge 3$ chosen independently to be in the random sequence with probability $1/\log n$ (this is also known as the Cram{\' e}r model).  One obvious difference is that all primes larger than $2$ must be odd, whereas a random sequence would surely contain many even numbers.  But if we could account for divisibility by small primes (such as $2$ in our example), would a modified random model describe accurately the behavior of prime numbers?

There are many ways in which we could try to make this theme precise.  For instance, the Riemann hypothesis predicts that $|\pi(x)-\text{li}(x)|$ is bounded by $C(\epsilon) x^{\frac 12+ \epsilon}$ for any $\epsilon >0$ and some constant $C(\epsilon)$.  Fluctuations of size about $\sqrt{x}$ are indeed what one would expect if we select random sets of integers with $n\ge 3$ included in the set with probability $1/\log n$.  Thus the Riemann hypothesis is, at a crude level, consistent with a random model of primes, although if we inspect the error term $\pi(x) -\text{li}(x)$ in finer detail then the influence of zeros of $\zeta(s)$ would be visible, and such features would deviate (in small but significant ways) from the random model.  

At the 1912 ICM, Landau posed four ``unattackable" problems on primes: (i) the Goldbach problem that every even integer larger than $2$ is the sum of two primes, (ii) the twin prime problem that there are infinitely many prime pairs $n$ and $n+2$, (iii) there is always a prime between two consecutive squares, and (iv) there are infinitely many primes of the form $n^2+1$.  All four problems remain open today, and all statements are exactly what one would expect for random sequences.  For example, the Cram{\' e}r model would suggest that the chance that $n$ and $n+2$ are both ``prime" is about $1/\log n \times 1/\log (n+2)$, which would predict about $x/(\log x)^2$ twin primes up to $x$.  Of course some care is needed, since the same prediction could be made for $n$ and $n+1$ being prime, and we will address this soon.  Similarly, we may expect that an even number $N$ may have about $N/(\log N)^2$ representations as a sum of two primes, making the Goldbach conjecture very plausible, and related arguments suggest the last two Landau problems as well.   

For the third Landau problem on the number of primes between $n^2$ and $(n+1)^2$, the random model already predicts what we believe to be the right answer --- namely, there should be about $(2n+1)/\log (n^2) \approx n/\log n$ primes in this interval.  For the other three problems, some modification must be made to the Cram{\' e}r model, to take into account the deterministic features of these problems with respect to divisibility by small primes.   Precise conjectures for these problems were first made by Hardy and Littlewood motivated by their work on the circle method.  These conjectures are widely believed to be true, and are supported by extensive heuristic and numerical evidence.   For instance, Hardy and Littlewood formulated the following conjecture for the number of twin primes below $x$: 
$$ 
\#\{ n \le x: \ \ n, n+2 \text{ both prime } \} \sim {\mathfrak S}(\{0, 2\}) \int_2^x \frac{dt}{(\log t)^2}. 
$$ 
Here $\int_2^x dt/(\log t)^2$ is asymptotically $x/(\log x)^2$, and corresponds to the prediction of the Cram{\' e}r model, while ${\mathfrak S}(\{0, 2\})$, known as the \emph{singular series}, is a correction factor 
$$ 
{\mathfrak S}(\{ 0, 2\}) = 2 \prod_{p\ge 3} \Big( 1-\frac 2p\Big) \Big(1 -\frac 1p \Big)^{-2} = 1.32 \ldots . 
$$ 
The constant ${\mathfrak S}(\{ 0, 2\})$ has a compelling probabilistic interpretation: it is a product over all primes $p$ (the first factor $2$ corresponds to the prime $p=2$), with the factor at $p$ keeping track of the ratio between the chance that $n$ and $n+2$ are not divisible by $p$, and the chance that two random numbers are not divisible by $p$.  Thus, for $p=2$, the chance that $n$ and $n+2$ are both not divisible by $2$ is $(1-1/2)$ ($n$ must be odd), while the chance that two random numbers are both not divisible by $2$ is $(1-1/2)^2 =1/4$; the ratio of these chances gives the correction factor $2$.   For primes $p\ge 3$, the chance that $n$ and $n+2$ are both not divisible by $p$ is $(1-2/p)$ whereas the chance that two random numbers are both not divisible by $p$ is $(1-1/p)^2$, and we see the corresponding correction factor in the definition of ${\mathfrak S}(\{ 0, 2\})$.  

Similar conjectures can be made for the binary Goldbach problem, or for the number of primes of the form $n^2+1$, modifying and correcting the naive predictions of the Cram{\' e}r model.    
To illustrate, we give a generalization of the conjecture for twin primes for counting prime $k$-tuples:  given 
distinct integers $h_1$, $h_2$, $\ldots$, $h_k$, for large $x$ how many integers $n \le x$ are there with $n+h_1$, $\ldots$, $n+h_k$ all being prime.  Here the Hardy--Littlewood conjecture predicts that 
\begin{equation} 
\label{1} 
\# \{ n\le x: \ \ n+h_1, \ldots, n+h_k \text{ all prime } \} \sim {\mathfrak S}(\{ h_1, \ldots , h_k \} ) \int_2^x \frac{dt}{(\log t)^k} 
\end{equation} 
where, with ${\mathcal H}= \{ h_1, \ldots, h_k\}$, 
\begin{equation} 
\label{2} 
{\mathfrak S}({\mathcal H}) = \prod_p \Big(1 -\frac{\nu({\mathcal H},p)}p\Big) \Big(1- \frac{1}{p}\Big)^{-k},
\end{equation} 
and $\nu({\mathcal H},p)$ denotes the number of distinct residue classes occupied by the set ${\mathcal H}$ viewed $\bmod \ p$.  Since $\nu({\mathcal H},p) =k$ if $p$ is larger than $\max |h_i-h_j|$, the product defining ${\mathfrak S}({\mathcal H})$ converges absolutely to a non-negative real number, and it equals zero only if $\nu({\mathcal H},p) =p$ 
for some prime $p$.  If $\nu({\mathcal H},p) =p$, then for any integer $n$ at least one of the numbers $n+h_1, \ldots, n+h_k$ would be a multiple of $p$, and therefore there can be only finitely many integers $n$ with $n+h_1$, $\ldots$, $n+h_k$ all being prime; for example this is what happens if we ask for $n$ and $n+1$ to be prime, or $n$, $n+2$, $n+4$ all to be prime.  When there is no such divisibility obstruction to $n+h_1$, $\ldots$, $n+h_k$ all being prime, the Hardy--Littlewood conjecture predicts a rich supply of such prime $k$-tuples.  This is perhaps the most central question in prime number theory, and remains open in any situation where ${\mathfrak S}({\mathcal H})$ is non-zero. 

\smallskip 

\noindent \textbf{Sieve theory.}  We have described quickly some of the main motivating questions in the theory of primes.  One main source of progress towards these questions is \emph{sieve theory}, and a large part of Maynard's work lies broadly in this area.   A typical problem in sieve theory is to bound the size of sets of integers ${\mathcal A}$ whose elements are constrained to omit $\nu(p)$ given residue classes $\bmod\ p$ for primes $p$.  For instance the twin prime problem is of this form, as we seek to find integers $n$ that are neither $0$ nor $-2\bmod \ p$ for all primes $p \le \sqrt{n+2}$ (so that $n$ and $n+2$ would both be prime).  In great generality sieve methods can produce upper bounds of the conjectured order of magnitude; for example, one can show that the number of twin primes up to $x$ is no more that $4$ times the conjectured Hardy--Littlewood asymptotic.  Producing corresponding lower bounds has proved to be a much harder problem, but sieve methods have led to striking partial results such as Chen's theorem that there are many primes $p$ for which $p+2$ has at most two prime factors, or Iwaniec's theorem that there are many $n$ for which $n^2+1$ has at most two prime factors.  For a comprehensive treatment of the subject, see \cite{FI2}.

Chen's theorem and Iwaniec's theorem exhibit a limitation of traditional sieve methods, known as the \emph{parity problem}, which often prevents us from knowing the parity of elements left unsieved, and thus from producing prime numbers.  But in some special cases, sieve methods in conjunction with other analytic input have produced prime numbers.  For instance, for large $x$ Baker, Harman, and Pintz \cite{BHP} showed that the interval $[x, x+x^{\theta}]$ contains at least $c x^{\theta}/\log x$ primes, where $c>0$ is a constant and $\theta = 0.525$;  the Landau problem of producing primes between consecutive squares corresponds to intervals with $\theta = \frac 12$.   Another spectacular example is due to Friedlander and Iwaniec \cite{FI1} who established an asymptotic formula for the number of primes up to $x$ that may be written as $n^2 +m^4$, an approximation to the Landau problem of primes of the form $n^2+1$.  A closely related result of Heath-Brown and Li \cite{HBLi} produces an asymptotic formula for primes of the form $n^2 +p^4$, where $p$ is prime.   Yet another beautiful result due to Heath-Brown \cite{HB1} establishes an asymptotic formula for the number of primes below $x$ of the form $n^3+2m^3$ with $m, n \in {\mathbb N}$.    Heath-Brown's result may be viewed as an approximation to the problem of producing primes of the form $n^3+2$, but before his work it was not even known if there are infinitely many primes that are the sum of three cubes of natural numbers!   A crucial feature of these results is that they deal with primes represented by specializations of \emph{norm forms}.  The Friedlander--Iwaniec result is concerned with the norm form $x^2+y^2 = N(x+iy)$ associated to the field ${\mathbb Q}(i)$, and specializing $y$ to be a square; Heath-Brown's result is concerned with the norm form $N(x+ y \alpha + z \alpha^2)$ taking the norm over the field ${\mathbb Q}(\alpha)$ with $\alpha= 2^{\frac 13}$, and specializing $z$ to be $0$.    The results of Friedlander--Iwaniec and Heath-Brown gave the first examples of thin sequences (in the sense that the number of integers below $X$ in the sequence is $\le X^{1-\delta}$ for some $\delta >0$) of polynomial values in two or more variables that represent infinitely many primes; no example is known of a  polynomial in $1$ variable of degree more than $1$ that represents infinitely many primes.  
 
 Maynard's work \cite{M4} gives a substantial generalization of Heath-Brown's approach, and produces many further examples of thin sequences of polynomial values in many variables that represent primes.  Consider an algebraic number $\omega \in {\mathbb C}$ of degree $n$, and let $K$ denote the field ${\mathbb Q}(\omega)$.  We can associate to this the norm form $N(\sum_{i=1}^{n} x_i \omega^{i-1})$, which is a homogeneous polynomial of degree $n$ in the variables $x_1$, $\ldots$, $x_n$.  A thin polynomial in many variables would be obtained by specializing some of the variables in this norm form to be zero; say, we set $x_{n-k+1}$, $\ldots$, $x_n=0$, and the number integers below $x$ represented by such an incomplete norm form would be about $x^{1-k/n}$.  In the range $n\ge 4k$, Maynard establishes an asymptotic formula for the number of primes represented by such an incomplete norm form, when the variables $x_1$, $\ldots$, $x_{n-k}$ take integer values in the range $[1, X]$. 
 
\smallskip 
 
\noindent \textbf{The circle method.}  Apart from sieve theory, another important source of progress towards problems on primes is the \emph{circle method}, which as we already mentioned formed the original motivation for Hardy and Littlewood in formulating their conjectures.  To illustrate, consider the Goldbach problem of representing an even integer $N$ as a sum of two primes.  Using Fourier analysis, the number of such representations of $N$ may be written as 
\begin{equation} 
\label{3} 
r(N) = \int_0^1 S(\alpha)^2 e^{-2\pi i N \alpha} d\alpha, \qquad \text{ where } \qquad S(\alpha) =\sum_{p\le N} e^{2\pi i p\alpha}. 
\end{equation} 
The idea in the circle method is that generating functions such as $S(\alpha)$ above tend to be large near rational numbers with small denominator (the \emph{major arcs}) and small away from them (the \emph{minor arcs}).    

While the circle method has not been able to tackle the binary Goldbach problem or the problem of twin primes, it has been extremely effective in problems where there is a bit more freedom.  For instance, the ternary Goldbach problem asks to represent odd numbers as a sum of three primes, and there is one extra variable to play with here.   Vinogradov famously used the circle method to show that all large odd numbers are the sum of three primes, and Helfgott \cite{Helf} has extended this to show that all odd numbers larger than $5$ may be so represented.  Here we may mention an impressive result of Matom{\" a}ki, Maynard, and Shao \cite{MMS} which shows that large odd numbers $n$ may be expressed as $p_1+ p_2 +p_3$, where all three primes $p_i$ lie in a short interval $[n/3 - n^{\theta}, n/3+ n^{\theta}]$ for any $\theta > 11/20$.  We mentioned earlier the work of Baker, Harman and Pintz \cite{BHP} showing the existence of primes in short intervals 
$[x,x+x^{0.525}]$, and the work of \cite{MMS} is remarkable in solving the ternary Goldbach problem using primes in only slightly longer intervals.  
 
A second example of what it might mean to have an extra degree of freedom is the Green--Tao theorem that the primes contain arbitrarily long arithmetic progressions $n$, $n+d$, $\ldots$, $n+(k-1)d$.  The Hardy--Littlewood conjecture would predict a stronger ``one dimensional'' version of such a result with specified choices for the common difference $d$; for instance, there should be infinitely many $k$-tuples primes of the form $n$, $n + k!$, $n+ 2 \cdot k!$, $\ldots$, $n+ (k-1) \cdot k!$.  The work of Green, Tao, and Ziegler \cite{GT1, GT2, GTZ} may be thought of as a far-reaching generalization of the circle method, obtaining asymptotic formulae for the number of prime solutions to linear systems with at ``least two degrees of freedom.''

Maynard's beautiful result on primes with missing digits (Theorem \ref{thm2} stated above) is a rare occasion where the circle method can be used to solve a binary problem.  Let ${\mathcal M}$ denote the set of natural numbers with no $7$ in their decimal expansion (naturally one could omit any other digit instead of $7$).  The number of integers in ${\mathcal M}$ up to $N$ is about $N^{\log 9/\log 10}= N^{1-\delta}$ with $\delta = 0.046\ldots$, so that ${\mathcal M}$ is a thin set making the problem of finding primes in it a challenge.  Before Maynard's work, Dartyge and Mauduit \cite{DM1, DM2} had used sieve theory to show that ${\mathcal M}$ contains integers with at most two prime factors.  To count the number of primes in ${\mathcal M}$ up to $N$, we use Fourier analysis writing this as 
$$ 
\sum_{\substack{ p \le N \\ p\in {\mathcal M}} } 1 = \int_0^1 S(\alpha) M(-\alpha) d\alpha, 
$$ 
where $S(\alpha)$ is the exponential sum over primes defined in \eqref{3}, and $M(\alpha) = \sum_{ m\le N, m \in {\mathcal M}} e^{2\pi i \alpha}$ is the corresponding exponential sum over the set ${\mathcal M}$.  Usually such a binary problem is hopeless to attack via the circle method --- the reason being that even most optimistically we may only expect ``square-root cancellation'' in the exponential sums $S(\alpha)$ and $M(-\alpha)$ for generic $\alpha$, and even that would produce an integrand of size $N^{\frac 12} \times N^{\frac 12 (1-\delta)}$,  which is bigger than the expected main term of size about $N^{1-\delta}/\log N$.  A crucial feature in this problem is that the set ${\mathcal M}$ has a very convenient structure which results in the exponential sum $M(\alpha)$ often being unusually small.  For instance, Maynard shows that its $L^1$-norm satisfies 
$$ 
\int_0^1 |M(\alpha)| d\alpha \ll N^{0.32}, 
$$ 
with the key point being that the exponent $0.32$ is smaller even than $(1-\delta)/2$, which is the optimistic square-root cancellation that we mentioned.  Such estimates raise the hope of being able to attack Theorem \ref{thm2}, and the main idea can be seen transparently in Maynard's expository article \cite{M6}, where he proves an easier version of Theorem \ref{thm2} treating primes missing a digit in base $b$ with $b$ sufficiently large.   The set of integers up to $N$ missing a digit in base $b$ has size about $N^{\log (b-1)/\log b}$, and so the problem becomes easier as the base $b$ gets larger.  Getting the base down to $10$ turns out to be a fiendishly difficult problem, and is arguably more significant psychologically than for any mathematical reason.  Maynard \cite{M2} tackles this brilliantly by introducing a number of new ideas, including ideas from the geometry of numbers, different aspects of sieve theory, and comparisons with a Markov process.  We may expect that even in base $3$ there should be infinitely many primes with a given digit missing; in base $2$, the only digit that might be omitted is $0$, and we find the problem of whether there are infinitely many Mersenne primes, which lies beyond reasonable mathematics.  We close this discussion by pointing out two other beautiful results on the digits of prime numbers which have elements in common with Maynard's work:  namely, work of Mauduit and Rivat \cite{MauRiv} which shows (in particular) that the sum of the decimal digits of primes is equally likely to be odd or even, and work of Bourgain \cite{Bour} which allows one to specify a small proportion of the binary digits of primes.

\smallskip 

\noindent \textbf{Gaps between primes.}   We now turn to a discussion of Maynard's most spectacular result --- the sun amidst small stars --- namely, Theorem \ref{thm1} above on finding many primes in bounded intervals.   To describe the recent history of this problem, let us first discuss how primes are spaced typically.  The prime number theorem tells us that the $n$-th prime $p_n$ is about $n \log n$, so that the average spacing between two consecutive primes, $p_{n+1}-p_n$, is about $\log p_n$.  What is the distribution of the normalized spacings $(p_{n+1}-p_n)/\log p_n$?  The Cram{\' e}r random model for primes would predict that these normalized spacings should behave like a Poisson process, and that for any fixed interval $[\alpha, \beta] \in {\mathbb R}_{\ge 0}$ 
\begin{equation} 
\label{4} 
\lim_{N \to \infty} \frac 1N \# \Big \{ n\le N: \ \ \frac{p_{n+1}-p_n}{\log p_n} \in [\alpha, \beta] \Big\} = \int_\alpha^\beta e^{-t} dt = e^{-\alpha} - e^{-\beta}. 
\end{equation} 
Gallagher \cite{Gal} showed that this prediction is also implied by the more refined Hardy--Littlewood conjectures, the key point being that the singular series constants ${\mathfrak S}({\mathcal H})$ (see \eqref{2}) are approximately $1$ (matching the naive Cram{\' e}r model) on average over $k$-element sets ${\mathcal H}$.   

This conjecture on the normalized spacings between primes is wide open. Indeed if we denote by ${\mathcal L}$ the set of limit points of $(p_{n+1}-p_n)/\log p_n$, then even the qualitative statement that ${\mathcal L} = [0, \infty]$ (which follows at once from \eqref{4}) is currently unknown.  By creating long strings of composite numbers, Westzynthius established that ${\mathcal L}$ contains $\infty$, but for a long time no other limit point was known (although Erd{\H o}s and Ricci had established that ${\mathcal L}$ has positive Lebesgue measure).  Dramatic progress was made in the 2005 with the path-breaking work of Goldston, Pintz, and Y\i ld\i r\i m \cite{GPY}, who showed that for any $\epsilon >0$ there are infinitely many $n$ with $p_{n+1}-p_n \le \epsilon \log p_n$.  Thus there are small gaps between primes in comparison to the average, and $0$ is now known to be in ${\mathcal L}$.  Before the work of Goldston, Pintz, and Y\i ld\i r\i m, it was only known that the difference between consecutive primes became smaller than about $\frac 14$ of the average spacing, and their work opened the door to later advances including Maynard's Theorem \ref{thm1}.

Suppose $h_1$, $\ldots$, $h_k$ are distinct integers with ${\mathfrak S}(\{ h_1, \ldots, h_k\})  >0$; such tuples are called \emph{admissible}, and for example $\{ k!, 2 \cdot k!, \ldots, k \cdot k!\}$ is admissible.   The Hardy--Littlewood conjecture predicts that there are infinitely many $n$ with $n+h_1$, $\ldots$, $n+h_k$ all being prime.   Instead of wanting all $k$ of these numbers to be prime, what if we only ask for at least two of them to be prime?   This would already show that infinitely often there are bounded gaps between consecutive prime numbers.  Suppose we could find  non-negative weights $w(n)$  with the property that for large $x$ and each $j=1$, $\ldots$, $k$, 
\begin{equation} 
\label{5} 
\sum_{\substack{ x\le n\le 2x \\ n+h_j \text{ prime }} } w(n) > \frac 1k \sum_{x \le n\le 2x} w(n).
\end{equation} 
Then summing \eqref{5} over all $j=1$, $\ldots$, $k$ we would obtain 
\begin{equation} 
\label{6}
\sum_{x \le n \le 2x} \#\{ 1\le j\le k: n+h_j \text{ prime} \} \ w(n) > \sum_{x \le n \le 2x} w(n), 
\end{equation} 
from which it would follow that there must be some $n$ with at least $2$ primes among $n+h_1$, $\ldots$, $n+h_k$.  
Thinking of the weights as giving a probability measure on $x\le n\le 2x$, we may interpret \eqref{6} as saying that the expected number of primes among the $n+h_j$ is greater than $1$, so that there must be $n$ with at least $2$ primes in this $k$-tuple. 

The difficult problem is to construct weights satisfying \eqref{5}, and natural choices for such weights are suggested by sieve theory, in particular the theory of the Selberg sieve.  The standard choice of Selberg sieve weights (which are used to give an upper bound for the number of prime $k$-tuples $n+h_1$, $\ldots$, $n+h_k$) takes the shape 
$$ 
w(n) = \Big( \sum_{ \substack { d| (n+h_1) \cdots (n+h_k) \\ d\le R} }\mu(d) \Big( \frac{\log R/d}{\log R}\Big)^k \Big)^2.
$$
Clearly $w(n) \ge 0$ always.  Expanding out the sum, the right side of \eqref{5} (the sum over all $n \in [x,2x]$) may be evaluated asymptotically so long as $R^2\le x^{1-\epsilon}$.  The left side of \eqref{5} is more involved, and relies on understanding the distribution of primes in arithmetic progressions with the modulus of the progression going up to $R^2$.  The Bombieri--Vinogradov theorem permits such an understanding (at a level comparable to what the Generalized Riemann Hypothesis would give) so long as $R^2 \le x^{\frac 12-\epsilon}$, so that $R$ is now constrained to be $\le x^{\frac 14-\epsilon}$.  For this choice of weights, the expected number of primes among the $n+h_j$ turns out to be about  $(2k/(k+1)) \log R/\log x$, so that with $R\le x^{1/4-\epsilon}$ one only expects to find $\frac 12$ a prime in the $k$-tuple. 

Although the Selberg sieve weights described above had been optimized for upper bounds in the prime $k$-tuple problem, Goldston, Pintz, and Y\i ld\i r\i m made the surprising discovery that there are better choices of weights for optimizing the ratio of the sums in \eqref{5}.  They considered weights of the form 
$$ 
w(n)  =  \Big( \sum_{ \substack { d| (n+h_1) \cdots (n+h_k) \\ d\le R} }\mu(d) \Big( \frac{\log R/d}{\log R}\Big)^{k+\ell} \Big)^2, 
$$
for a suitable parameter $\ell$, which turns out in the optimal case to be around $\sqrt{k}$.  With this choice of weights, they found that the expected number of primes among $n+h_j$ is about twice as large as previously, being   $(4+O(1/k^{\frac 12})) \log R/\log x$.   With $R= x^{\frac 14-\epsilon}$, this barely fails to give the desired relation \eqref{5}, and thus barely falls short of proving bounded gaps between primes.  By considering an additional possible prime value $n+h$ for $1\le h \le \epsilon \log x$, Goldston, Pintz, Y\i ld\i r\i m were able to deduce from this argument that there are infinitely many $n$ with $p_{n+1} -p_n \le \epsilon \log n$.  For a more detailed discussion of these ideas see \cite{S}.

 If one could take $R$ to be $x^{\frac 14 +\delta}$ for any $\delta >0$, then the argument of Goldston, Pintz, and 
 Y\i ld\i r\i m would give bounded gaps between primes.  To take such a value for $R$, one would need to understand the distribution of primes up to $x$ in arithmetic progressions, when the modulus of the progression is as large as $x^{\frac 12+2\delta}$.  The Elliott--Halberstam conjectures predict that such results should hold (on average) when the modulus is as large as $x^{1-\epsilon}$.  Partial progress towards such extensions of the Bombieri--Vinogradov theorem was made by Fouvry and Iwaniec \cite{FoIw}, and Bombieri, Friedlander, and Iwaniec  \cite{BFI}, but these results did not apply immediately to the problem of showing bounded gaps between primes.   In April 2013, Yitang Zhang \cite{Zhang} made a spectacular breakthrough by establishing a version of the Bombieri--Vinogradov theorem in an extended range which was sufficient for the method of Goldston, Pintz, and Y\i ld\i r\i m.  Zhang established that if $k > 3.5 \times 10^6$ then for any admissible $k$-tuple $h_1$, $\ldots$, $h_k$ there are infinitely many $n$ with at least two of the $n+h_j$ being prime.  This implied that infinitely often the gaps between consecutive primes is less than $70$ million.   Refinements of Zhang's work on the equidistribution of primes in arithmetic progressions were made by the Polymath project \cite{Polymath}, and still further qualitative and quantitative refinements of such results may be found in the recent papers of Maynard \cite{May1, May2, May3}.

 Zhang's work established the case $m=2$ of Theorem \ref{thm1}.  However, even if one could take the largest possible range for $R$, namely $R= x^{\frac 12 -\epsilon}$ (which would be permitted by the Elliott--Halberstam conjecture), the Goldston--Pintz--Y\i ld\i r\i m weights would only yield that the expected number of primes in an admissible $k$-tuple is 
 $\ge 2 -\epsilon$.  In other words, even under the Elliott--Halberstam conjecture one would fall short of establishing the existence of three primes in bounded intervals. 
 
 The proof of Theorem \ref{thm1} is based on a different choice of the weights $w(n)$, discovered just months after Zhang's work by Maynard (who announced the results in a memorable talk at Oberwolfach in October 2013) and independently by Tao (in unpublished work).   The Maynard--Tao weights are a multi-dimensional extension of the weights considered earlier, and take (roughly speaking) the shape
  $$ 
  w(n) = \Big( \sum_{ \substack{ d_1, \ldots, d_k \\ d_i | n+h_i \\ \prod d_i \le R }} \prod_{i=1}^{k} \mu(d_i) F\Big( 
  \frac{\log d_1}{\log R}, \ldots, \frac{\log d_k}{\log R}\Big) \Big)^2, 
  $$ 
  for suitable smooth functions $F: [0,1]^k \to {\mathbb R}$.  Astonishingly it turns out that for an appropriate choice for $F$, the expected number of primes in the tuple $n+h_1$, $\ldots$, $n+h_k$ (recall \eqref{6} above) is $\ge c \log k \frac{\log R}{\log x}$, for a positive constant $c$; in fact $c$ may be taken close to $1$ if $k$ is large enough.  
 The key point is that this expected number of primes in $k$-tuples tends to infinity with $k$, and in fact we only need $R$ to grow like any power of $x$ for the method to succeed, so that Bombieri--Vinogradov which permits $R=x^{\frac 14-\epsilon}$ is already sufficient!  Thus the following more precise version of Theorem \ref{thm1} holds, which may be viewed as a partial result towards the Hardy--Littlewood prime $k$-tuples conjecture.
 
 \begin{theorem} [Maynard \cite{M1}]   Let $m \ge 2$ be a natural number.  Let $k$ be sufficiently large in terms of $m$, and let ${\mathcal H} = \{ h_1, \ldots, h_k\}$ be any set of $k$ integers with ${\mathfrak S}({\mathcal H}) > 0$.  Then there exist infinitely many $n$ such that the $k$-tuple $n+h_1$, $\ldots$, $n+ h_k$ contains at least $m$ primes. 
 \end{theorem} 
 
 Maynard showed that $k$ may be taken smaller than $C m^2 e^{4m}$ for a suitable constant $C$, and further refinements of this (incorporating also the work of Zhang) have been made in the work of Baker and Irving \cite{BI} 
 who showed that $k$ may be taken as $C e^{3.815 m}$.  Of special interest is the case $m=2$ where the Polymath project \cite{Polymath2} optimized these arguments to establish that any admissible $50$-tuple  contains $2$ primes infinitely often.  In particular, they showed that $p_{n+1}- p_n \le 246$ infinitely often, and conditional on the Elliott--Halberstam conjecture that infinitely often there are at least two primes in the triple $n$, $n+2$, $n+6$.  Let us mention one other uniform variant of these results:  Maynard \cite{M5} shows, for instance, that there are at least $c X\exp(-\sqrt{\log X})$ values of  $x \in [X,2X]$ such that the interval $[x, x +\log X]$ contains at least $c\log \log X$ primes (here $c$ is a positive constant).  
 For detailed expositions on these results of Zhang, Maynard, and Tao, see \cite{Gran, Kow}.

The Maynard--Tao weights offer a flexible new method to study many problems on primes and related sequences, and have found a number of applications.    We describe two other results using these weights, both still concerned with spacings between consecutive primes.    We referred earlier to the result of Westzynthius on large gaps between consecutive primes, which showed that $\infty$ lies in the set ${\mathcal L}$ of limit points of the normalized spacings $(p_{n+1}-p_n)/\log p_n$.  This was quantified in the 1930's by Erd{\H o}s and Rankin who showed that, for a positive constant $C$ 
\begin{equation} 
\label{7} 
\max_{p_n \le X} (p_{n+1} - p_n) \ge C \log X \frac{(\log \log X) \log \log \log \log X}{(\log \log \log X)^2}. 
\end{equation} 
The random model would suggest that the maximal gap between primes up to $X$ should be about $(\log X)^2$.  This is known as Cram{\' e}r's conjecture, and while this is very delicate, it is widely believed that the maximal gap is no more than $(\log X)^{2+\epsilon}$, although even this is far beyond Landau's unattackable problem of the existence of a prime between consecutive squares.  Erd{\H o}s drew attention to the problem of finding larger gaps between consecutive primes, offering \$ 10 000 for a bound that would replace $C$ in \eqref{7} with a function tending to $\infty$ with $X$.  For more than 75 years, this problem resisted attack, with only improvements of the constant $C$ being known.   Then, by a remarkable coincidence, in 2014 \emph{two} different techniques emerged, both establishing \eqref{7} with $C$ replaced by a function tending to infinity with $X$.   One approach, by Ford, Green, Konyagin, and Tao \cite{FGKT}, built upon the work of Green--Tao on arithmetic progressions in the primes, while the other approach, by Maynard \cite{M3}, found a way to adapt the Maynard--Tao sieve weights.  The second approach was better suited for quantifying the large gaps that are produced, and, joining forces, Ford, Green, Konyagin, Maynard, and Tao \cite{FGKMT} established that for some constant $C>0$ 
\begin{equation} 
\label{8} 
\max_{p_n \le X} (p_{n+1} - p_n) \ge C \log X \frac{(\log \log X) \log \log \log \log X}{\log \log \log X}, 
\end{equation} 
 improving the bound in \eqref{7} by a factor of $\log \log \log X$.

The results on small gaps and large gaps between consecutive primes show that $0$ and $\infty$ lie in the set ${\mathcal L}$ of limit points of the normalized prime spacings.  No other explicit numbers are known to lie in ${\mathcal L}$, although we expect ${\mathcal L}$ to include all non-negative real numbers.   Following Zhang's breakthrough, Pintz \cite{Pintz2} showed that ${\mathcal L}$ contains an interval $[0,c]$ for some $c>0$, which however is ineffective and cannot be computed explicitly.  Using the Maynard--Tao sieve weights, Banks, Freiberg, and Maynard \cite{BFM} established the following beautiful result:  If $\beta_1 \le \beta_2 \le \ldots \le \beta_9$ are any nine real numbers, then at least one of their differences $\beta_j -\beta_i$ (with $i<j$) must be an element of ${\mathcal L}$.    Their result has been refined by Pintz \cite{Pintz}, and Merikoski \cite{Mer}, and Merikoski shows that the same result holds if we start with just four real numbers $\beta_1 \le \beta_2 \le \beta_3\le \beta_4$.  Moreover, Merikoski has also shown that for any $T >0$, the set ${\mathcal L} \cap [0,T]$ has measure at least $T/3$.

\smallskip 

\noindent \textbf{The Duffin--Schaeffer conjecture.}  So far we have focussed entirely on Maynard's work concerned with prime numbers.  In a very different direction, Maynard in joint work with Koukoulopoulos \cite{KM}, resolved one of the central problems in the metric theory of Diophantine approximation, known as the Duffin--Schaeffer conjecture.  

Diophantine approximation is concerned with finding rational approximations $a/q$ to a given irrational number $\alpha$, with an emphasis on making $|\alpha -a/q|$ small in terms of $q$.  The most basic result is Dirichlet's theorem that for every irrational number $\alpha$, there are infinitely many rational approximations $a/q$, with $a\in {\mathbb Z}$, $q\in {\mathbb N}$ and $(a,q)=1$ (so that the fraction is in reduced form) such that $|\alpha -a/q| \le 1/q^2$.  For quadratic irrationals (like $\sqrt{2}$ or the golden ratio), Dirichlet's theorem is essentially the best possible, and for every such $\alpha$ there exists a positive constant $C(\alpha)$ such that 
$|\alpha -a/q| \ge C(\alpha)/q^2$ for any rational approximation $a/q$. A celebrated result of Roth establishes that for any algebraic irrational $\alpha$ and any $\epsilon >0$ one has $|\alpha -a/q| \ge C(\alpha, \epsilon)/q^{2+\epsilon}$, for a suitable positive constant $C(\alpha,\epsilon)$.  For particular interesting transcendental numbers, such as $\pi$, it remains an outstanding open problem to determine how well they can be approximated by rational numbers. 

Metric Diophantine approximation is concerned with such approximation problems that hold for \emph{almost all} irrational numbers $\alpha$, with \emph{almost all} interpreted in the sense of Lebesgue measure.  Since the problem of approximating $\alpha$ by rationals is identical to that of approximating $\alpha+1$, we may restrict attention to irrational numbers $\alpha \in [0,1)$.  The most basic problem is the following:  suppose $\psi: {\mathbb N} \to {\mathbb R}_{\ge 0}$ is a given function, what can be said about the measure of $\alpha \in [0,1)$ for which there exist infinitely many rational numbers $a/q$ in reduced form (that is, $(a,q)=1$) with $|\alpha- a/q| \le \psi(q)$.  For instance, Dirichlet's theorem tells us that if $\psi(q)=1/q^2$, then all irrational $\alpha \in [0,1)$ admit infinitely many such rational approximations.  

Let ${\mathcal A}_q = {\mathcal A}_q(\psi)$ denote the set of $\alpha\in [0,1)$ for which there exists some reduced fraction $a/q$ with $|\alpha- a/q| \le \psi(q)$, and let ${\mathcal A}$ denote the set of $\alpha \in [0,1)$ lying in infinitely many of the sets ${\mathcal A}_q$.  Thus 
$$ 
{\mathcal A} = \bigcap_{Q=1}^{\infty} {\widetilde {\mathcal A}}(Q), \qquad \text{ with } \qquad {\widetilde {\mathcal A}}(Q) = \bigcup_{q=Q}^{\infty} {\mathcal A}_q. 
$$ 
Now the measure of ${\mathcal A}_q$ is $\le 2 \phi(q) \psi(q)$, since there are $\phi(q)$ possible choices for the numerator $a$, and if $\psi(q) \le 1/(2q)$ so that the intervals for different $a$ do not overlap then equality holds here.  If 
$\sum_{q=1}^{\infty} \phi(q) \psi(q)$ converges, then the measure of $\widetilde{\mathcal A}(Q)$ is bounded by $2\sum_{q=Q}^{\infty} \phi(q) \psi(q)$, which is the tail of a convergent series and thus tends to $0$ as $Q\to \infty$.  
It follows that ${\mathcal A}$ has measure $0$.  This argument is identical to the easy part of the Borel--Cantelli Lemma.

In 1941, Duffin and Schaeffer made the remarkable conjecture that in the complementary case when $\sum_{q=1}^{\infty} \phi(q) \psi(q)$ diverges, the measure of ${\mathcal A}$ is $1$.  Since then the Duffin--Schaeffer conjecture has remained one of the central motivating questions in the theory of metric Diophantine approximations.  A number of partial results towards this conjecture were established: for example, a beautiful result of Gallagher \cite{Gal2} showed that the measure of the set ${\mathcal A}(\psi)$ is always either $0$ or $1$, work of Erd{\H o}s \cite{E} and Vaaler \cite{V} established the conjecture when $\psi(q)$ is $O(1/q^2)$ for all $q$, higher dimensional analogues of the conjecture were proved by Pollington and Vaughan \cite{PV},  and weaker versions of the conjecture with extra divergence conditions were established in \cite{Harman, HPV, Ais}.   But the full problem resisted until the recent work of Koukoulopoulos and Maynard \cite{KM}: 

\begin{theorem} [Koukoulopoulos and Maynard \cite{KM}]  Let $\psi: {\mathbb N} \to {\mathbb R}_{\ge 0}$ be such that 
$\sum_{q=1}^{\infty} \phi(q) \psi(q)$ diverges.  Then the set of $\alpha \in [0,1)$ that have infinitely many rational approximations $|\alpha -a/q| \le \psi(q)$ with $(a,q)=1$ has Lebesgue measure $1$.  In other words, the Duffin--Schaeffer conjecture holds.  
\end{theorem} 

We refer to Koukoulopoulos's talk at this ICM \cite{K} for a more detailed exposition of this result, and the ideas behind its proof. 

\smallskip 

We have given an overview of some of Maynard's most spectacular achievements in analytic number theory.  Maynard's work is characterized by ingenious but simple ideas, which are carried very far with his powerful technical ability.  As impressive as his work so far has been, it may only mark a beginning.



\begin{funding}
This work was partially supported by grants from the National Science Foundation, and a Simons Investigator Award from the Simons Foundation.  \end{funding}


\bibliographystyle{emss}

\bibliography{LaudRefs}{}








\end{document}